\documentclass[reqno]{amsart}
\usepackage{amsfonts}
\usepackage{amsmath}
\usepackage{amssymb}
\usepackage{hyperref}
\usepackage[square,numbers,sort&compress]{natbib}

\newtheorem{theorem}{Theorem}[section]

\theoremstyle{definition}
\newtheorem{definition}[theorem]{Definition}
\newtheorem{remark}[theorem]{Remark}

\AtBeginDocument{{\noindent\small \sc{}\newline
}\vspace{9mm}}
\input{tcilatex}

\begin{document}
\title[On the modified $q$-Genocchi numbers and polynomials]{On the modified 
$q$-Genocchi numbers and polynomials and their applications\\
\hspace{0.5cm}}
\author{Serkan Araci, Armen Bagdasaryan, Erkan A\u{g}y\"{u}z and Mehmet
Acikgoz}
\address{(S. Araci)\\
Atat\"{u}rk Street, 31290 Hatay, Turkey }
\email{mtsrkn@hotmail.com}
\address{(A. Bagdasaryan)\\
Russian Academy of Sciences, Institute for Control Sciences\\
65 Profsoyuznaya, 117997 Moscow, Russia }
\email{abagdasari@hotmail.com}
\address{(E. A\u{g}y\"{u}z and M. Acikgoz)\\
University of Gaziantep, Faculty of Arts and Science, Department of
Mathematics, 27310 Gaziantep, Turkey. }
\email{erkanagyuz@hotmail.com}
\email{acikgoz@gantep.edu.tr}
\date{}

\begin{abstract}
The main objective of this paper is to introduce the modified $q$-Genocchi 
polynomials and to define their generating function. In the paper, we show
new relations, which are explicit formula, derivative formula,
multiplication formula, and some others, for mentioned $q$-Genocchi
polynomials. By applying Mellin transformation to the generating function of
the modified $q$-Genocchi polynomials, we define $q$-Genocchi zeta-type
functions which are interpolated by the modified $q$-Genocchi polynomials at
negative integers.
\end{abstract}

\date{November 23, 2013}
\maketitle


\pagenumbering{arabic}


\section{Introduction, Definitions and Notations}

As is well known, Genocchi polynomials are given by%
\begin{equation}
\frac{2t}{e^{t}+1}e^{xt}=\sum_{n=0}^{\infty }G_{n}\left( x\right) \frac{t^{n}%
}{n!}\text{, }\left\vert t\right\vert <\pi \text{.}  \label{eq.1}
\end{equation}

Putting $x=0$ into (\ref{eq.1}), we have the quantity $G_{n}\left( 0\right)
:=G_{n}$ that stands for Genocchi numbers, some values of which are as
follows: 
\begin{equation*}
G_{0}=0,G_{1}=1,G_{2}=-1,G_{4}=1,G_{6}=-3,G_{8}=17,G_{10}=-155,G_{12}=2073,%
\cdots 
\end{equation*}%
with $G_{2n+1}=0$ for $n\geq 1$. The even coefficients of Genocchi numbers
can be computed by%
\begin{equation*}
G_{2n}=2\left( 1-2^{2n}\right) B_{2n}
\end{equation*}%
where $B_{n}$ implies Bernoulli numbers and the last identity of the above
is known as \textit{Genocchi's }theorem (see \cite{Araci}, \cite{Araci1}, 
\cite{Cangul}, \cite{Horadam}, \cite{Horadam 1}, \cite{Srivastava}). One of
the most recent papers on the theory of Genocchi numbers and polynomials is
the paper of A. F. Horadam \cite{Horadam}, which deals mainly with the
theory of Genocchi polynomials. Information on Bernoulli polynomials $%
B_{n}\left( x\right) $ and Euler polynomials $E_{n}\left( x\right) $, to
which $G_{n}\left( x\right) $ may be related, has been derived in [12-14,
18-22, 24, 25]. While a lot of the properties of Genocchi polynomials bear a
striking resemblance to the properties of Bernoulli and Euler polynomials,
some properties are rather different. Note that Genocchi polynomials occur
naturally in the areas of elementary number theory, complex analytic number
theory, homotopy theory (stable homotopy groups of spheres), differential
topology (differential structures on spheres), theory of modular forms
(Eisenstein series), $p$-adic analytic number theory ($p$-adic $L$%
-functions), quantum physics (quantum groups).

We need the following notations and definitions for the sequel of this paper:

Suppose that $p$ be a fixed odd prime number. Throughout this paper, we will
employ the following notations:%
\begin{eqnarray*}
&&%
\mathbb{Z}
_{p},\text{ the ring of }p\text{-adic rational integers} \\
&&%
\mathbb{Q}
_{p},\text{ the field of }p\text{-adic rational numbers} \\
&&%
\mathbb{C}
_{p},\text{ denotes the completion of algebraic closure of }%
\mathbb{Q}
_{p} \\
&&%
\mathbb{N}
,\text{ the set of natural numbers.}
\end{eqnarray*}

In addition to the notation above, we will use the quantity $%
\mathbb{N}
^{\ast }$ meaning $%
\mathbb{N}
\cup \left\{ 0\right\} $.

Let $\upsilon _{p}$ be normalized exponential valuation of $%
\mathbb{C}
_{p}$ such that 
\begin{equation*}
\left\vert p\right\vert _{p}=p^{-\upsilon _{p}\left( p\right) }=\frac{1}{p}%
\text{.}
\end{equation*}

When one talks of $q$-extension, $q$ will be known as an indeterminate,
either a complex number $q\in 
\mathbb{C}
$, or a $p$-adic number $q\in 
\mathbb{C}
_{p};$ it is always clear from the context. If $q\in 
\mathbb{C}
$, then one usually assumes that $\left\vert q\right\vert <1$. If $q\in 
\mathbb{C}
_{p}$, then one usually assumes that $\left\vert q-1\right\vert _{p}<1,$ and
hence $q^{x}=\exp \left( x\log q\right) $ for $x\in 
\mathbb{Z}
_{p}$ (see \cite{Araci1}, \cite{Araci2}, \cite{Araci3}, \cite{Cangul}, \cite%
{Seo}, \cite{Kim2}, \cite{Kim5}, \cite{Kim6}, \cite{D. S. Kim1}, \cite{Lee}, 
\cite{Rim}, \cite{Rim1}, \cite{Simsek1}, \cite{Simsek2}).

In the theory of $q$-calculus for a real parameter $q\in \left( 0,1\right) $%
, $q$-numbers are given by 
\begin{equation*}
\left[ x\right] _{q}=\frac{1-q^{x}}{1-q}\text{ \textup{and} }\left[ x\right]
_{-q}=\frac{1-\left( -q\right) ^{x}}{1+q}\text{,}
\end{equation*}%
(for details, see \cite{Kac}). Note that $\lim_{q\rightarrow 1}\left[ x%
\right] _{q}=x$ for any $x$ with $\left\vert x\right\vert _{p}\leq 1$ in the
present $p$-adic case.

Let $UD\left( 
\mathbb{Z}
_{p}\right) $ be the space of uniformly differentiable functions on $%
\mathbb{Z}
_{p}$. For a positive integer $d$ with $\left( d,p\right) =1$, set 
\begin{eqnarray*}
X &=&X_{d}=\lim_{\overleftarrow{n}}%
\mathbb{Z}
/dp^{n}%
\mathbb{Z}
\text{, }X_{1}=%
\mathbb{Z}
_{p} \\
X^{\ast } &=&\underset{\underset{\left( a,p\right) =1}{0<a<dp}}{\cup }\left(
a+dp%
\mathbb{Z}
_{p}\right) 
\end{eqnarray*}%
and%
\begin{equation*}
a+dp^{n}%
\mathbb{Z}
_{p}=\left\{ x\in X\mid x\equiv a\left( \func{mod}dp^{n}\right) \right\} 
\text{,}
\end{equation*}%
where $a\in 
\mathbb{Z}
$ satisfies the condition $0\leq a<dp^{n}$ cf. \cite{Araci1}, \cite{Araci2}, 
\cite{Araci3}, \cite{Cangul}, \cite{Seo}, \cite{Kim2}, \cite{Kim5}, \cite%
{Kim6}, \cite{D. S. Kim1}, \cite{Lee}, \cite{Rim}, \cite{Rim1}, \cite%
{Simsek1}, \cite{Simsek2}.

The following $p$-adic $q$-Haar distribution is defined by Kim \cite{Kim6} as%
\begin{equation*}
\mu _{q}(x+p^{n}%
\mathbb{Z}
_{p})=\frac{q^{x}}{[p^{n}]_{q}}\text{.}
\end{equation*}

Thus, for $f\in UD\left( 
\mathbb{Z}
_{p}\right) $, the $p$-adic $q$-integral on $%
\mathbb{Z}
_{p}$ was defined by Kim as follows:%
\begin{eqnarray}
I_{q}\left( f\right) &=&\int_{%
\mathbb{Z}
_{p}}f\left( x\right) d\mu _{q}\left( x\right)  \label{equation 2} \\
&=&\lim_{n\rightarrow \infty }\frac{1}{\left[ p^{n}\right] _{q}}%
\sum_{x=0}^{p^{n}-1}f\left( x\right) q^{x}\text{.}  \notag
\end{eqnarray}

On the one hand, the bosonic integral is considered as the bosonic limit $%
q\rightarrow 1,$ $I_{1}\left( f\right) =\lim_{q\rightarrow 1}I_{q}\left(
f\right) $, which are called \textit{Volkenborn} integral. On the other
hand, the \textit{fermionic} $p$-adic $q$-integral on $%
\mathbb{Z}
_{p}$ \ is considered by Kim as follows:%
\begin{equation}
I_{-q}\left( f\right) =\lim_{t\rightarrow -q}I_{t}\left( f\right) =\int_{%
\mathbb{Z}
_{p}}f\left( x\right) d\mu _{-q}\left( x\right) \text{.}  \label{equation 3}
\end{equation}

By (\ref{equation 3}), we have the following well-known integral equation 
\begin{equation}
qI_{-q}\left( f_{1}\right) +I_{-q}\left( f\right) =\left[ 2\right]
_{q}f\left( 0\right) \text{,}  \label{equation 4}
\end{equation}%
where $f_{1}\left( x\right) $ is a translation with $f\left( x+1\right) $
cf. \cite{Araci3}, \cite{Kim1}.

Recently, the modified $q$-Bernoulli polynomials are introduced by \cite{Seo}%
\begin{equation*}
\widetilde{\beta }_{n,q}\left( x\right) =\int_{%
\mathbb{Z}
_{p}}\left( x+\left[ y\right] _{q}\right) ^{n}d\mu _{q}\left( y\right) \text{
for }n\in 
\mathbb{N}
^{\ast }\text{.}
\end{equation*}

In \cite{Rim1}, the modified $q$-Euler polynomials are given by%
\begin{equation*}
\widetilde{\epsilon }_{n,q}\left( x\right) =\int_{%
\mathbb{Z}
_{p}}q^{-y}\left( x+\left[ y\right] _{q}\right) ^{n}d\mu _{-q}\left(
y\right) \text{ for }n\in 
\mathbb{N}
^{\ast }\text{.}
\end{equation*}

These polynomials have interesting properties for doing study in the theory
of Riemann-zeta function and Euler-zeta function. By the same motivation,
let us now consider the following definition.

\begin{definition}
For $n\in 
\mathbb{N}
^{\ast }$, the modified $q$-Genocchi polynomials are defined by 
\begin{equation}
\frac{\mathcal{G}_{n+1,q}\left( x\right) }{n+1}=\int_{%
\mathbb{Z}
_{p}}\left( x+\left[ y\right] _{q}\right) ^{n}d\mu _{-q}\left( y\right) 
\text{ for }n\in 
\mathbb{N}
^{\ast }\text{.}  \label{equation 5}
\end{equation}

\begin{remark}
Setting $q\rightarrow 1$ in (\ref{equation 5}), it yields 
\begin{equation*}
\lim_{q\rightarrow 1}\left( \frac{\mathcal{G}_{n+1,q}\left( x\right) }{n+1}%
\right) :=\frac{G_{n+1}\left( x\right) }{n+1}=\int_{%
\mathbb{Z}
_{p}}\left( x+y\right) ^{n}d\mu _{-1}\left( y\right) \text{ }.
\end{equation*}
\end{remark}
\end{definition}

In the next section, we derive some new interesting identities for the
modified $q$-Genocchi polynomials and get the generating function for the
modified $q$-Genocchi polynomials. Next, by applying Mellin transformation
to the generating function of the modified $q$-Genocchi polynomials, we
define an analogue of Hurwitz Genocchi zeta function. It shows that this
zeta function is interpolated by the modified $q$-Genocchi polynomials at
negative integers.

\section{On the modified $q$-Genocchi numbers and polynomials}

In this part, we give the properties of the modified $q$-Genocchi numbers
and polynomials and define the generating function of the modified $q$%
-Genocchi polynomials.

Substituting $x=0$ in (\ref{equation 5}), we have%
\begin{equation}
\frac{\mathcal{G}_{n+1,q}\left( 0\right) }{n+1}:=\frac{\mathcal{G}_{n+1,q}}{%
n+1}=\int_{%
\mathbb{Z}
_{p}}\left[ y\right] _{q}^{n}d\mu _{-q}\left( y\right) \text{ for }n\in 
\mathbb{N}
^{\ast }\text{.}  \label{equation 6}
\end{equation}

It is clear that $\mathcal{G}_{0,q}=0$. By (\ref{equation 5}) and (\ref%
{equation 6}), we have the following theorem.

\begin{theorem}
(\textbf{Addition formula}) For $n\in 
\mathbb{N}
^{\ast }$, we get%
\begin{equation*}
\mathcal{G}_{n,q}\left( x\right) =\sum_{k=0}^{n}\binom{n}{k}x^{n-k}\mathcal{G%
}_{k,q}=\left( \mathcal{G}_{q}+x\right) ^{n}
\end{equation*}%
where we have used the technique of umbral calculus as $\left( \mathcal{G}%
_{q}\right) ^{n}:=$ $\mathcal{G}_{n,q}$.
\end{theorem}

On account of (\ref{equation 5}), we have%
\begin{eqnarray*}
\frac{\mathcal{G}_{n+1,q}\left( x\right) }{n+1} &=&\int_{%
\mathbb{Z}
_{p}}\left( x+\left[ y\right] _{q}\right) ^{n}d\mu _{-q}\left( y\right)  \\
&=&\left( \frac{1}{1-q}\right) ^{n}\int_{%
\mathbb{Z}
_{p}}\left( \left( 1-q\right) x+1-q^{y}\right) ^{n}d\mu _{-q}\left( y\right) 
\\
&=&\sum_{k=0}^{n}\binom{n}{k}\left( \frac{1}{1-q}\right) ^{k}x^{n-k}\int_{%
\mathbb{Z}
_{p}}\left( 1-q^{y}\right) ^{k}d\mu _{-q}\left( y\right)  \\
&=&\sum_{k=0}^{n}\sum_{j=0}^{k}\binom{n}{k}\binom{k}{j}\left( \frac{1}{1-q}%
\right) ^{k}x^{n-k}\left( -1\right) ^{j}\int_{%
\mathbb{Z}
_{p}}q^{jy}d\mu _{-q}\left( y\right)  \\
&=&\left[ 2\right] _{q}\sum_{k=0}^{n}\sum_{j=0}^{k}\binom{n}{k}\binom{k}{j}%
\left( \frac{1}{1-q}\right) ^{k}\frac{x^{n-k}\left( -1\right) ^{j}}{1+q^{j+1}%
}\text{.}
\end{eqnarray*}

Therefore, we arrive at the following theorem that establishes an explicit
formula for the modified $q$-Genocchi polynomials.

\begin{theorem}
(\textbf{Explicit formula}) For $n\in 
\mathbb{N}
^{\ast }$, we have%
\begin{equation}
\frac{\mathcal{G}_{n+1,q}\left( x\right) }{n+1}=\left[ 2\right]
_{q}\sum_{k=0}^{n}\sum_{j=0}^{k}\binom{n}{k}\binom{k}{j}\left( \frac{1}{1-q}%
\right) ^{k}\frac{x^{n-k}\left( -1\right) ^{j}}{1+q^{j+1}}\text{.}
\label{equation 8}
\end{equation}
\end{theorem}

Let $F_{q}\left( x;t\right) =\sum_{n=0}^{\infty }\mathcal{G}_{n,q}\left(
x\right) \frac{t^{n}}{n!}$ where $F_{q}\left( x;t\right) $ can be written by
(\ref{equation 5}), as follows:%
\begin{equation*}
t\sum_{n=0}^{\infty }\left( \int_{%
\mathbb{Z}
_{p}}\left( x+\left[ y\right] _{q}\right) ^{n}d\mu _{-q}\left( y\right)
\right) \frac{t^{n}}{n!}=t\int_{%
\mathbb{Z}
_{p}}e^{\left( x+\left[ y\right] _{q}\right) t}d\mu _{-q}\left( y\right) 
\text{.}
\end{equation*}

That is, the generating function of the modified $q$-Genocchi polynomials
can be rewritten as%
\begin{equation}
\sum_{n=0}^{\infty }\mathcal{G}_{n,q}\left( x\right) \frac{t^{n}}{n!}=t\int_{%
\mathbb{Z}
_{p}}e^{\left( x+\left[ y\right] _{q}\right) t}d\mu _{-q}\left( y\right) 
\text{.}  \label{equation 7}
\end{equation}

On the other hand, by (\ref{equation 8}), we see that%
\begin{eqnarray*}
\sum_{n=0}^{\infty }\mathcal{G}_{n,q}\left( x\right) \frac{t^{n}}{n!} &=&%
\left[ 2\right] _{q}t\sum_{n=0}^{\infty }\left( \sum_{k=0}^{n}\sum_{j=0}^{k}%
\binom{n}{k}\binom{k}{j}\left( \frac{1}{1-q}\right) ^{k}\frac{x^{n-k}\left(
-1\right) ^{j}}{1+q^{j+1}}\right) \frac{t^{n}}{n!} \\
&=&\left[ 2\right] _{q}te^{xt}\left( \sum_{n=0}^{\infty }\left( \frac{1}{1-q}%
\right) ^{n}\sum_{j=0}^{n}\binom{n}{j}\frac{\left( -1\right) ^{j}}{1+q^{j+1}}%
\frac{t^{n}}{n!}\right) \\
&=&\left[ 2\right] _{q}t\sum_{m=0}^{\infty }\left( -1\right)
^{m}q^{m}e^{\left( x+\left[ m\right] _{q}\right) t}\text{.}
\end{eqnarray*}

Thus, we state the following theorem.

\begin{theorem}
\label{thm gen.}(\textbf{Generating function}) We have two forms of
representation for the generating function of the modified q-Genocchi
polynomials as follows:%
\begin{eqnarray*}
F_{q}\left( x;t\right) &=&\sum_{n=0}^{\infty }\mathcal{G}_{n,q}\left(
x\right) \frac{t^{n}}{n!} \\
&=&t\int_{%
\mathbb{Z}
_{p}}e^{\left( x+\left[ y\right] _{q}\right) t}d\mu _{-q}\left( y\right) \\
&=&\left[ 2\right] _{q}t\sum_{m=0}^{\infty }\left( -1\right)
^{m}q^{m}e^{\left( x+\left[ m\right] _{q}\right) t}\text{.}
\end{eqnarray*}
\end{theorem}

By Theorem \ref{thm gen.}, we arrive at%
\begin{equation*}
\sum_{n=0}^{\infty }\mathcal{G}_{n,q}\left( x\right) \frac{t^{n}}{n!}%
=\sum_{n=0}^{\infty }\left[ \left[ 2\right] _{q}\sum_{m=0}^{\infty }\left(
-1\right) ^{m}q^{m}\left( x+\left[ m\right] _{q}\right) ^{n}\right] \frac{%
t^{n+1}}{n!}\text{.}
\end{equation*}

Equating the coefficients of $t^{n}$ in the last identity above, we get the
following theorem.

\begin{theorem}
For $n\in 
\mathbb{N}
^{\ast }$, we have%
\begin{equation}
\frac{\mathcal{G}_{n+1,q}\left( x\right) }{n+1}=\left[ 2\right]
_{q}\sum_{m=0}^{\infty }\left( -1\right) ^{m}q^{m}\left( x+\left[ m\right]
_{q}\right) ^{n}\text{.}  \label{equation 9}
\end{equation}
\end{theorem}

Let us take $\frac{d}{dx}$, which is a derivative operator, in the both
sides of expression in Theorem \ref{thm gen.}; we see that 
\begin{equation*}
\frac{d}{dx}\left[ \mathcal{G}_{n,q}\left( x\right) \right] =n\mathcal{G}%
_{n-1,q}\left( x\right) \text{.}
\end{equation*}

Therefore we obtain the following theorem.

\begin{theorem}
\label{thm deriv.}(\textbf{Derivative property for} $\mathcal{G}_{n,q}\left(
x\right) $) For $n\in 
\mathbb{N}
$, we have%
\begin{equation*}
\frac{d}{dx}\left[ \mathcal{G}_{n,q}\left( x\right) \right] =n\mathcal{G}%
_{n-1,q}\left( x\right) \text{.}
\end{equation*}
\end{theorem}

Polynomials $A_{n}\left( x\right) $ are called \textit{Appell polynomials} 
\cite{Fort}, \cite{D. S. Kim}, \cite{D. S. Kim1} if they have the property 
\begin{equation*}
\frac{dA_{n}\left( x\right) }{dx}=nA_{n-1}\left( x\right) \text{.}
\end{equation*}

Thus, by Theorem \ref{thm deriv.}, the aforementioned modified $q$-Genocchi
polynomials are \textit{Appell polynomials}.

Let $d$ be an odd positive integer. We are now in a position to state the
property of multiplication for the modified $q$-Genocchi polynomials with
the help of the aforementioned fermionic $p$-adic $q$-integral on $%
\mathbb{Z}
_{p}$ in this paper:%
\begin{eqnarray*}
\frac{\mathcal{G}_{n+1,q}\left( x\right) }{n+1} &=&\int_{%
\mathbb{Z}
_{p}}\left( x+\left[ y\right] _{q}\right) ^{n}d\mu _{-q}\left( y\right)  \\
&=&\lim_{m\rightarrow \infty }\frac{1}{\left[ dp^{m}\right] _{-q}}%
\sum_{y=0}^{dp^{m}-1}\left( x+\left[ y\right] _{q}\right) ^{n}\left(
-q\right) ^{y} \\
&=&\lim_{m\rightarrow \infty }\frac{1}{\left[ dp^{m}\right] _{-q}}%
\sum_{y=0}^{p^{m}-1}\sum_{a=0}^{d-1}\left( x+\left[ a+dy\right] _{q}\right)
^{n}\left( -q\right) ^{a+dy} \\
&=&\frac{1}{\left[ d\right] _{-q}}\sum_{a=0}^{d-1}\left( -q\right)
^{a}\lim_{m\rightarrow \infty }\frac{1}{\left[ p^{m}\right] _{-q^{d}}}%
\sum_{y=0}^{p^{m}-1}\left( x+\left[ a\right] _{q}+q^{a}\left[ d\right] _{q}%
\left[ y\right] _{q^{d}}\right) ^{n}\left( -q\right) ^{dy} \\
&=&\frac{\left[ d\right] _{q}^{n}}{\left[ d\right] _{-q}}\sum_{a=0}^{d-1}%
\left( -1\right) ^{a}q^{a\left( n+1\right) }\left[ \lim_{m\rightarrow \infty
}\frac{1}{\left[ p^{m}\right] _{-\left( q^{d}\right) }}\sum_{y=0}^{p^{m}-1}%
\left( \frac{x+\left[ a\right] _{q}}{q^{a}\left[ d\right] _{q}}+\left[ y%
\right] _{q^{d}}\right) ^{n}\left( -q^{d}\right) ^{y}\right]  \\
&=&\frac{\left[ d\right] _{q}^{n}}{\left[ d\right] _{-q}}\sum_{a=0}^{d-1}%
\left( -1\right) ^{a}q^{a\left( n+1\right) }\int_{%
\mathbb{Z}
_{p}}\left( \frac{x+\left[ a\right] _{q}}{q^{a}\left[ d\right] _{q}}+\left[ y%
\right] _{q^{d}}\right) ^{n}d\mu _{-q^{d}}\left( y\right) 
\end{eqnarray*}%
where we have used the following two identities that are well known in the
theory of $q$-calculus 
\begin{equation*}
\left[ x+y\right] _{q}=\left[ x\right] _{q}+q^{x}\left[ y\right] _{q}\text{
and }\left[ xy\right] _{q}=\left[ x\right] _{q}\left[ y\right] _{q^{x}}\text{
(see \cite{Kac}).}
\end{equation*}

As a result of these applications, we have the following theorem.

\begin{theorem}
(\textbf{Multiplication formula}) Let $d\equiv 1\left( \func{mod}2\right) $
and $n\in 
\mathbb{N}
$, then we have%
\begin{equation*}
\mathcal{G}_{n,q}\left( q^{a}\left[ d\right] _{q}x\right) =\frac{\left[ d%
\right] _{q}^{n-1}}{\left[ d\right] _{-q}}\sum_{a=0}^{d-1}\left( -1\right)
^{a}q^{a\left( n+1\right) }\mathcal{G}_{n,q^{d}}\left( x+\frac{\left[ a%
\right] _{q}}{q^{a}\left[ d\right] _{q}}\right) \text{.}
\end{equation*}
\end{theorem}

\section{On the modified $q$-Genocchi polynomials in connection with
Zeta-type function}

The Hurwitz zeta function is defined by 
\begin{equation}
\zeta \left( s,x\right) =\sum_{n=0}^{\infty }\frac{1}{\left( n+x\right) ^{s}}%
\text{, for }s\in 
\mathbb{C}
\text{.}  \label{equation 10}
\end{equation}

Putting $x=1$ in (\ref{equation 10}), yields to%
\begin{equation*}
\zeta \left( s,1\right) :=\zeta \left( s\right) =\sum_{n=1}^{\infty }\frac{1%
}{n^{s}}
\end{equation*}%
where $\zeta \left( s\right) $ is known as Riemann zeta function to be
convergence for $\func{Re}\left( s\right) >1$ (cf. \cite{Bagdasaryan1}, \cite%
{Bagdasaryan2}, \cite{Bagdasaryan3}, \cite{Sri}, \cite{Kim4}). Note that the
Bernoulli numbers interpolate by the Riemann zeta function, which plays a
crucial role in analytic number theory and has many applications in physics,
probability and applied statistics. Firstly, Leonard Euler introduced the
Riemann zeta function of real argument without using complex analysis. By (%
\ref{equation 10}), we have the following relation: For $n\in 
\mathbb{N}
$,%
\begin{equation*}
\zeta \left( 1-n\right) =-\frac{B_{n}}{n}\text{ (see \cite{Bagdasaryan1}, 
\cite{Bagdasaryan2}, \cite{Bagdasaryan3}, \cite{Sri}, \cite{Kim4}).}
\end{equation*}

So we consider an analogue of zeta function by applying Mellin
transformation to the generating function of the modified $q$-Genocchi
polynomials. From those consideration, we write the following:%
\begin{equation*}
\frac{1}{\Gamma \left( s\right) }\int_{0}^{\infty }t^{s-2}F_{q}\left(
x;-t\right) dt
\end{equation*}%
this identity yields to%
\begin{equation}
\left[ 2\right] _{q}\sum_{m=0}^{\infty }\frac{\left( -1\right) ^{m}q^{m}}{%
\left( x+\left[ m\right] _{q}\right) ^{s}}\text{.}  \label{equation 11}
\end{equation}

As $q$ approaches to $1$ in (\ref{equation 11}), it becomes%
\begin{equation}
2\sum_{m=0}^{\infty }\frac{\left( -1\right) ^{m}}{\left( x+m\right) ^{s}}%
\text{.}  \label{equation 12}
\end{equation}

The Equation (\ref{equation 12}) is known as Genocchi zeta function \cite%
{Simsek1}. It shows that the Equation (\ref{equation 11}) seems to be the
new $q$-analogue of Genocchi zeta function. So we give the following
definition.

\begin{definition}
Let $s\in 
\mathbb{C}
$ with $\func{Re}\left( s\right) >1$ and $0<x\leq 1$. The new $q$-analogue
of the Hurwitz Genocchi zeta-type function is expressed by%
\begin{equation}
\mathcal{\zeta }_{q}\left( s,x\right) =\left[ 2\right] _{q}\sum_{m=0}^{%
\infty }\frac{\left( -1\right) ^{m}q^{m}}{\left( x+\left[ m\right]
_{q}\right) ^{s}}\text{.}  \label{equation 13}
\end{equation}
\end{definition}

By (\ref{equation 9}) and (\ref{equation 13}), we have 
\begin{equation}
\mathcal{\zeta }_{q}\left( -n,x\right) =\frac{\mathcal{G}_{n+1,q}\left(
x\right) }{n+1}\text{.}  \label{equation 14}
\end{equation}

The Equation (\ref{equation 14}) seems to be interpolating for the modified $%
q$-Genocchi polynomials at negative integers. Now we define partial $q$%
-Hurwitz Genocchi zeta-type function as follows: For $F\equiv 1\left( \func{%
mod}2\right) $ 
\begin{equation}
\mathcal{H}_{q}\left( s,x:a,F\right) =\left[ 2\right] _{q}\sum_{\underset{m>0%
}{m\equiv a\left( \func{mod}F\right) }}\frac{\left( -1\right) ^{m}q^{m}}{%
\left( x+\left[ m\right] _{q}\right) ^{s}}\text{.}  \label{equation 15}
\end{equation}

By (\ref{equation 15}), we have%
\begin{eqnarray*}
\mathcal{H}_{q}\left( s,x:a,F\right) &=&\left[ 2\right] _{q}\sum_{m=0}^{%
\infty }\frac{\left( -1\right) ^{mF+a}q^{mF+a}}{\left( x+\left[ mF+a\right]
_{q}\right) ^{s}} \\
&=&\left( -1\right) ^{a}q^{a}\left[ 2\right] _{q}\sum_{m=0}^{\infty }\frac{%
\left( -1\right) ^{m}q^{mF}}{\left( x+\left[ a\right] _{q}+q^{a}\left[ F%
\right] _{q}\left[ m\right] _{q^{F}}\right) ^{s}} \\
&=&\frac{\left( -1\right) ^{a}q^{a\left( 1-s\right) }\left[ 2\right] _{q}}{%
\left[ F\right] _{q}^{s}\left[ 2\right] _{q^{F}}}\left[ \left[ 2\right]
_{q^{F}}\sum_{m=0}^{\infty }\frac{\left( -1\right) ^{m}q^{mF}}{\left( \frac{%
x+\left[ a\right] _{q}}{q^{a}\left[ F\right] _{q}}+\left[ m\right]
_{q^{F}}\right) ^{s}}\right] \\
&=&\frac{\left( -1\right) ^{a}q^{a\left( 1-s\right) }\left[ 2\right] _{q}}{%
\left[ F\right] _{q}^{s}\left[ 2\right] _{q^{F}}}\mathcal{\zeta }%
_{q^{F}}\left( s,\frac{x+\left[ a\right] _{q}}{q^{a}\left[ F\right] _{q}}%
\right) \text{.}
\end{eqnarray*}

Therefore we procure the following theorem.

\begin{theorem}
For $F\equiv 1\left( \func{mod}2\right) $ and $0\leq a<F$, we have 
\begin{equation}
\mathcal{H}_{q}\left( s,x:a,F\right) =\frac{\left( -1\right) ^{a}q^{a\left(
1-s\right) }\left[ 2\right] _{q}}{\left[ F\right] _{q}^{s}\left[ 2\right]
_{q^{F}}}\mathcal{\zeta }_{q^{F}}\left( s,\frac{x+\left[ a\right] _{q}}{q^{a}%
\left[ F\right] _{q}}\right) \text{.}  \label{equation 16}
\end{equation}
\end{theorem}

Putting $s=-n$ in (\ref{equation 16}) yields to 
\begin{equation*}
\mathcal{H}_{q}\left( -n,x:a,F\right) =\frac{\left( -1\right) ^{a}q^{a\left(
n+1\right) }\left[ 2\right] _{q}\left[ F\right] _{q}^{n}}{\left[ 2\right]
_{q^{F}}}\frac{\mathcal{G}_{n+1,q^{F}}\left( \frac{x+\left[ a\right] _{q}}{%
q^{a}\left[ F\right] _{q}}\right) }{n+1}\text{.}
\end{equation*}

\end{document}